\documentclass[12pt]{article}

\newtheorem{theorem}{Theorem}[section]
\newtheorem{lemma}[theorem]{Lemma}
\newtheorem{proposition}[theorem]{Proposition}
\newtheorem{corollary}[theorem]{Corollary}

\newenvironment{poof}{\textit{Proof:  }}{
~\hfill\rule{2mm}{3mm}\vspace{.2in}}

\usepackage{amssymb,amsmath,tabularx,graphicx,float,placeins}
\usepackage{tikz}
\usetikzlibrary[calc,intersections,through,backgrounds,arrows,decorations.pathmorphing]

\def\Ical{\mathcal{I}}

\def\Nbb{\mathbb{N}}

\def\Rbb{\mathbb{R}}
\def\Sbb{\mathbb{S}}
\def\Zbb{\mathbb{Z}}

\def\ds{\displaystyle}
\def\ra{\rightarrow}

\def\ov{\overline}

\def\pr{\prime}

\def\setm{\setminus}

\DeclareMathOperator{\Asc}{Asc}

\DeclareMathOperator{\Int}{Int}
\DeclareMathOperator{\susp}{susp}
\DeclareMathOperator{\height}{ht}

\DeclareMathOperator{\nonc}{nonc}
\DeclareMathOperator{\prop}{prop}

\begin{document}

\title{Homotopy type of intervals of the second higher Bruhat orders}
\date{\today}
\author{Thomas McConville}

\maketitle

\begin{abstract}
The higher Bruhat order is a poset of cubical tilings of a cyclic zonotope whose covering relations are cubical flips. For a 2-dimensional zonotope, the higher Bruhat order is a poset on commutation classes of reduced words for the longest element of a type A Coxeter system. For this case, we prove that the noncontractible intervals are in natural correspondence with the zonogonal tilings of a zonogon. Our proof uses some tools developed by Felsner and Weil to show that the two standard orderings on the rhombic tilings of a zonogon are identical.
\end{abstract}

\section{Introduction}\label{sec_introduction}

\begin{figure}
\centering
\includegraphics{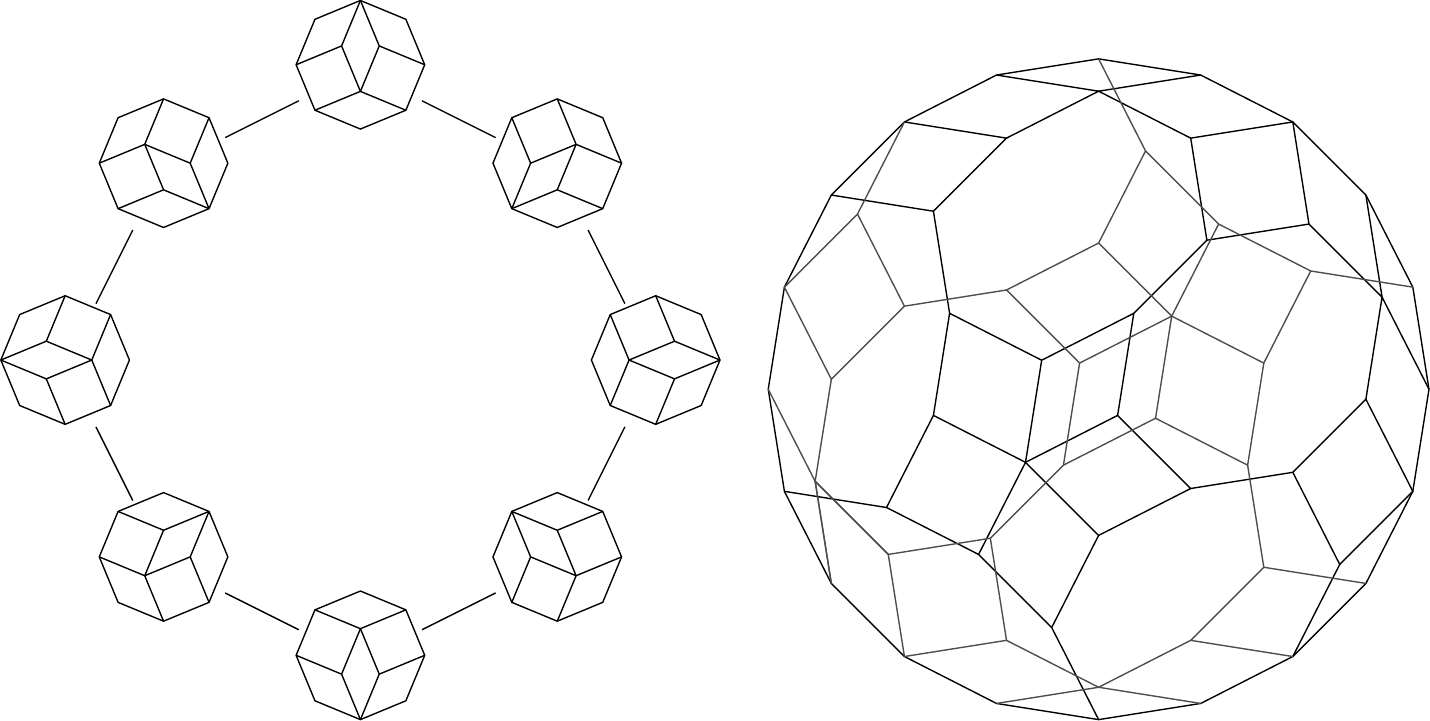}
\caption{\label{fig_B42}(left) $B(4,2)$ as a poset of rhombic tilings of a zonogon.  (right) $B(5,2)$}
\end{figure}

The higher Bruhat order is a poset structure on the cubical tilings of a cyclic zonotope, ordered by upward flips (Figure \ref{fig_B42}).  Alternatively, the higher Bruhat order $B(n,d)$ is the poset of consistent subsets, defined as follows.

Let $\binom{[n]}{d+1}$ denote the $(d+1)$-element subsets of $\{1,\ldots,n\}$.  We say a subset $X$ of $\binom{[n]}{d+1}$ is \emph{closed} if $I\cup\{i,j\}\in X$ and $I\cup\{j,k\}\in X$ implies $I\cup\{i,k\}\in X$ for $I\in\binom{[n]}{d-1},\ i,j,k\in[n]-I,\ i<j<k$.  For instance, $\{123,134\}$ is not a closed subset of $\binom{[4]}{3}$ since it contains $\{1\}\cup\{2,3\}$ and $\{1\}\cup\{3,4\}$ but not $\{1\}\cup\{2,4\}$.  A subset $X$ of $\binom{[n]}{d+1}$ is \emph{consistent} (or \emph{biclosed} or \emph{clopen}) if both $X$ and $\binom{[n]}{d+1}-X$ are closed.  The consistent sets are ordered by \emph{single-step inclusion}; that is $X\leq Y$ holds if there exists a sequence $X=X_0\subseteq\cdots\subseteq X_t=Y$ of consistent sets for which $|X_i-X_{i-1}|=1$ for all $i$.  In particular, $B(n,1)$ may be identified with the weak order on the symmetric group on $[n]$.  The second higher Bruhat order $B(n,2)$ defines an ordering on the commutation-equivalence classes of reduced words for the longest element of the symmetric group on $[n]$.


Introduced by Manin and Schechtman, the higher Bruhat orders have many equivalent interpretations, including single-element extensions of an alternating matroid, cubical tilings of a cyclic zonotope, and ``admissible'' permutations of $d$-subsets of $[n]$ up to a suitable equivalence; see \cite[Theorem 4.1]{ziegler:bruhat}, \cite{manin.schechtman:arrangements}, or \cite{kapranov.voevodsky:categories}.  The higher Bruhat orders have appeared in a wide variety of contexts, including higher categories and Zamolodchikov's tetrahedral equation \cite{kapranov.voevodsky:categories}, soliton solutions of the Kadomtsev-Petviashvili equation \cite{dimakis.muller-hoissen:kp}, and the multidimensional cube recurrence \cite{henriques.speyer:multidimensional}.

We consider the homotopy type of intervals of $B(n,2)$.  Our main result is that every interval of $B(n,2)$ is either contractible or homotopy equivalent to a sphere.  As usual, the topology associated to a poset $P$ is that of its order complex, the simplicial complex of chains $x_0<\cdots<x_m$ of elements of $P$.  If $P$ is a bounded poset, we usually consider the order complex of its \emph{proper part} $P_{\prop}$, the same poset with those bounds removed.  As the M\"obius invariant of $P$ is equal to the reduced Euler characteristic of the order complex of $P_{\prop}$, the homotopy types of intervals of a poset completely determines the M\"obius function.  This is a useful point of view, as many techniques for determining M\"obius functions have homotopy analogues \cite{bjorner:topological}.

\begin{figure}
\centering
\includegraphics{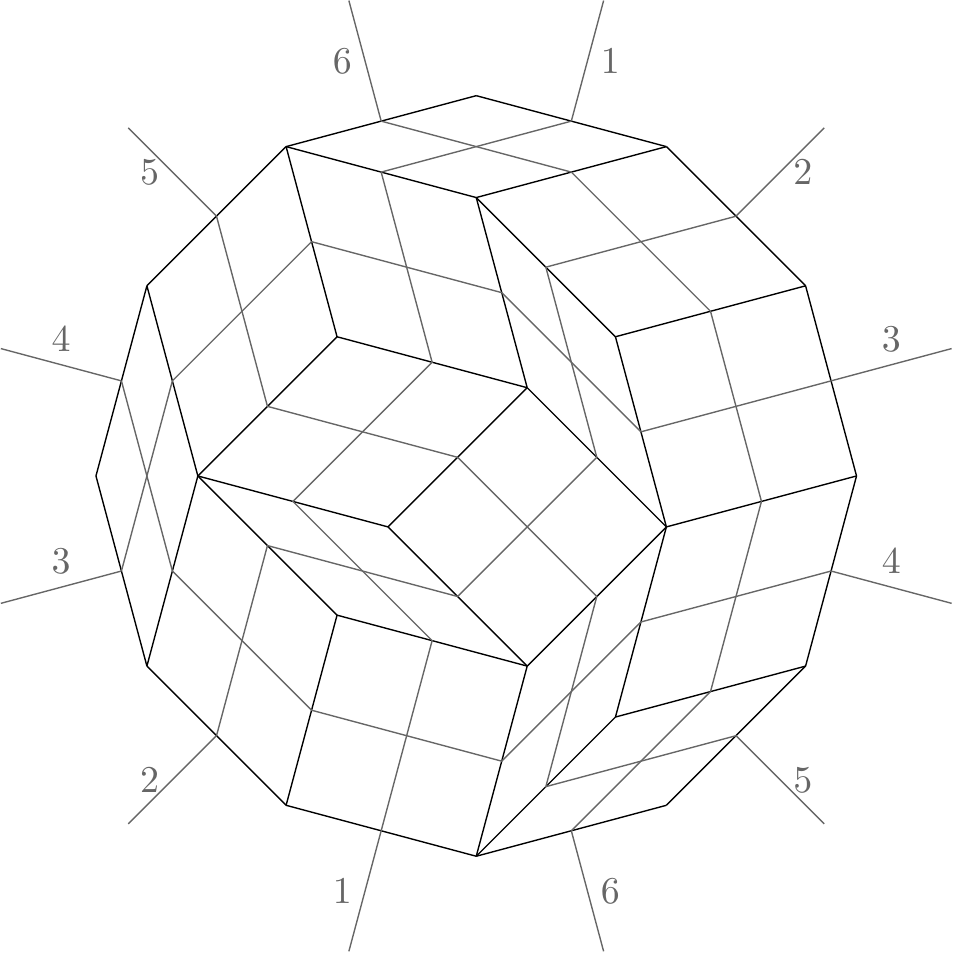}
\caption{\label{fig_bohne_dress}A rhombic tiling of a zonogon with its associated pseudoline arrangement.}
\end{figure}

Consistent subsets of $\binom{[n]}{3}$ are in natural bijection with simple pseudoline arrangements with $n$ pseudolines, cyclically ordered at infinity.  We also identify one of the two infinite regions bounded by $1$ and $n$ as the ``bottom'' region.  The consistent set associated to a simple pseudoline arrangement is the set of inversions of the arrangement, where $\{i<j<k\}\in\binom{[n]}{3}$ is an \emph{inversion} if the crossing of the pseudolines $i$ and $k$ occurs below $j$.  The inversion set of the pseudoline arrangement in Figure \ref{fig_bohne_dress} is $\{124,134,135,234,235\}$.  Simple pseudoline arrangements also correspond to rhombic tilings of a zonogon via the Bohne-Dress Theorem (\cite{bohne:kombinatorische},\cite{richter.ziegler:zonotopal}) as demonstrated in Figure \ref{fig_bohne_dress}.

\begin{figure}
\centering
\includegraphics{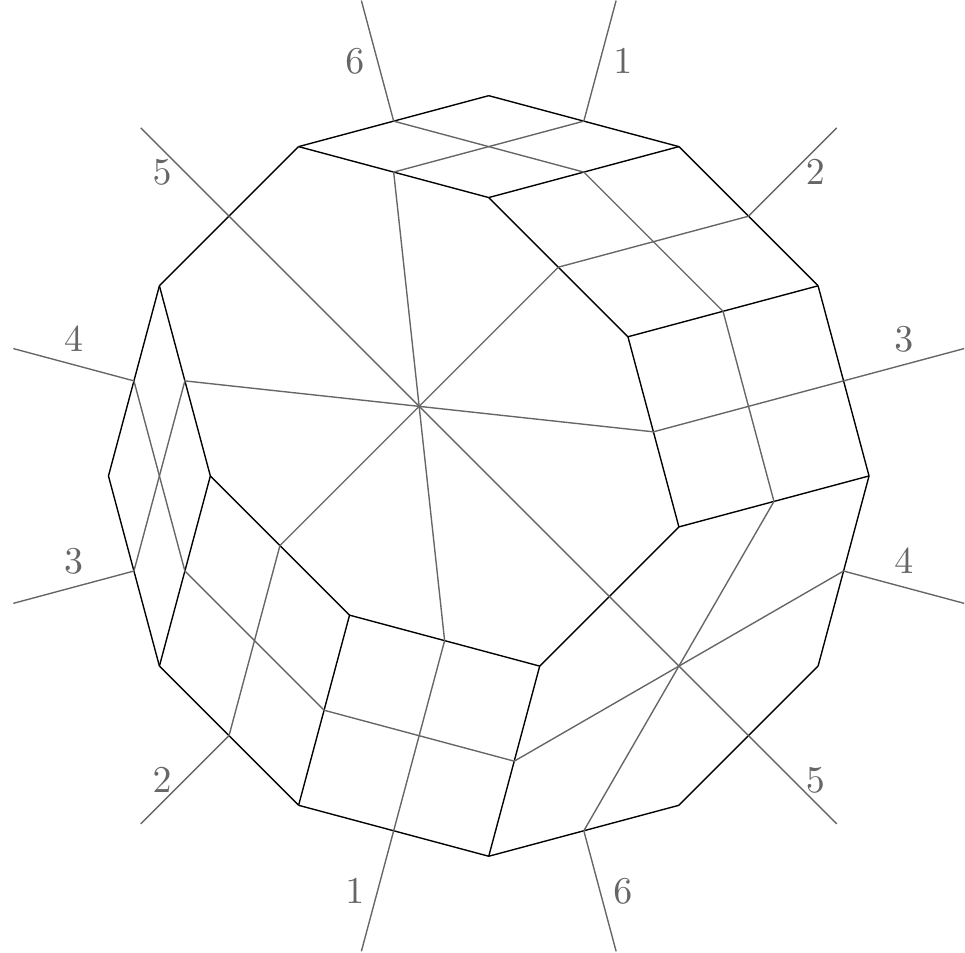}
\caption{\label{fig_zono_tile}A zonogonal tiling with its associated non-simple pseudoline arrangement.}
\end{figure}

A (non-simple) arrangement of pseudolines may have crossings involving more than two pseudolines.  The set of simple arrangements that may be obtained by resolving these crossings forms a closed interval of $B(n,2)$, which we call a \emph{facial interval}.  For example, the arrangement in Figure \ref{fig_zono_tile} has two non-simple crossings that may be resolved in 16 ways, which is an interval of $B(6,2)$.  One such resolution is the arrangement of Figure \ref{fig_bohne_dress}.


Rambau proved that that the proper part of $B(n,d)$ is homotopy equivalent to an $(n-d-2)$-sphere as an application of his Suspension Lemma \cite{rambau:suspension}.  Reiner extended this by showing that any facial interval of $B(n,d)$ is homotopy equivalent to a sphere \cite[Conjecture 6.9]{reiner:baues}.  He conjectured that every other interval is contractible.  As $B(n,1)$ is isomorphic to the weak order of the symmetric group on $[n]$, the conjectured homotopy type of intervals was already verified by Bj\"orner for $B(n,1)$ \cite{bjorner:orderings}.  We prove Reiner's conjecture for $B(n,2)$.

\begin{theorem}\label{thm_main}
An interval of $B(n,2)$ is non-contractible if and only if it is facial.
\end{theorem}

Bj\"orner's computation of the homotopy type of intervals of $B(n,1)$ relies on the lattice property of the weak order.  Indeed, the weak order is a crosscut-simplicial lattice, which means for any interval $(x,y)$, the join of any proper subset of atoms of $(x,y)$ is not equal to $y$.  Hence, every interval is either contractible or homotopy equivalent to a sphere by the Crosscut Lemma.

Although $B(n,2)$ is not a lattice when $n\geq 6$, it is ``close enough'' to being a lattice that a similar argument may be applied.  For any poset $P$, the order complex of $P$ is homotopy equivalent to the order complex of $P_{\nonc}$, the subposet of elements $X$ for which $\{Y\in P:\ Y<X\}$ is non-contractible \cite[Proposition 6.1]{walker:homotopy}.  We prove that if $P$ is any open interval of $B(n,2)$, then either $P_{\nonc}$ is the proper part of a Boolean lattice, or $P_{\nonc}$ contains an element $X$ such that $X\vee Y$ exists in $P_{\nonc}$ for all $Y\in P_{\nonc}$.  The latter intervals are contractible by a join-contraction argument.



By Theorem \ref{thm_main}, there is a poset isomorphism between the non-contractible intervals of $B(n,2)$ ordered by inclusion and the lifting space of a central arrangement of $n$ lines, ordered by weak maps.  Using some general techniques in poset topology, this isomorphism implies that the lifting space is homotopy equivalent to a sphere of dimension $n-3$.  In general, the lifting space of an alternating matroid is known to be homotopy equivalent to a sphere (\cite{sturmfels.ziegler:extension} Theorem 4.12); see Figure \ref{fig_B42}.  Reiner's conjecture would provide an alternate proof of this result.


The paper is organized as follows. Some topological preliminaries are given in Section \ref{sec_poset}. We prove some results on general higher Bruhat orders in Section \ref{sec_bruhat}. Wiring diagrams are defined in Section \ref{sec_bruhat_2} along with other results specific to the second higher Bruhat orders. Finally, the proof of Theorem \ref{thm_main} is given in Section \ref{sec_main}.

\section{Poset Topology}\label{sec_poset}


We establish some notation and recall a few fundamental results on the topology of posets, following Bj\"orner \cite{bjorner:topological}.  A \emph{poset} $P$ is a set with a reflexive, antisymmetric, transitive binary relation. We will tacitly assume that a given poset is finite unless specified otherwise. A \emph{lower bound} (\emph{upper bound}) is an element $\hat{0}$ ($\hat{1}$) such that $\hat{0}\leq x$ ($x\leq\hat{1}$) for all $x\in P$.  We say $P$ is \emph{bounded} if it has both an upper and lower bound.  A poset is a \emph{lattice} if every pair of elements $x,y\in P$ has a least upper bound $x\vee y$ and greatest lower bound $x\wedge y$.

If $P$ has an upper or a lower bound, then the \emph{proper part} $\ov{P}$ is the same poset with those bounds removed.  Given $x\leq y$, the \emph{closed interval} $[x,y]$ (\emph{open interval} $(x,y)$) is the set of $z\in P$ such that $x\leq z\leq y$ ($x<z<y$).

The \emph{order complex} $\Delta(P)$ of a poset $P$ is the abstract simplicial complex with vertex set $P$ and simplices $\{x_0,\ldots,x_d\}$ where $x_0<\cdots<x_d$ is a chain of $P$.  We define the topology of a poset to be that of its order complex.

The \emph{M\"obius function} $\mu:\Int(P)\ra\Zbb$ is the unique function on the closed intervals of $P$ such that
$$\mu([x,y])=\begin{cases}1\hspace{5mm}&\mbox{if }x=y\\\ds\sum_{z\in[x,y]}\mu([x,z])=0\ &\mbox{if }x<y\end{cases}$$
The value of $\mu([x,y])$ is the reduced Euler characteristic of $\Delta((x,y))$.  Many techniques for computing M\"obius functions have homotopy counterparts.

\begin{lemma}[\cite{walker:homotopy} Proposition 6.1]\label{lem_contractible_link}
If $x\in P$ such that $P_{<x}$ or $P_{>x}$ is contractible, then $P$ is homotopy equivalent to $P-\{x\}$.
\end{lemma}

Let $x_1<\cdots<x_N$ be a linear extension of $P$.  By deleting elements $x_i$ such that $P_{>x_i}$ is contractible in the order of the linear extension, we deduce the following result.

\begin{lemma}\label{lem_noncontractible}
$P$ is homotopy equivalent to the subposet
$$\{x\in P:\ P_{>x} \mbox{ is non-contractible } \}$$.
\end{lemma}

We let $P_{\nonc}$ be the subposet of Lemma \ref{lem_noncontractible}. Replacing $P_{>x}$ with $P_{<x}$, a dual form of this lemma also holds, but we will not need it.

Let $\Int(P)$ be the poset of closed intervals of $P$, ordered by inclusion.  It is known that $\ov{\Int}(P)$ is homeomorphic to the suspension of $\ov{P}$. This was originally proved by Walker \cite[Theorem 6.1(c)]{walker:homeomorphisms} by specifying a ``subdivision map'' between geometric realizations of their order complexes. An alternative proof was given in \cite[Lemma 3.3.10]{mcconville:thesis} by constructing the order complex of $\ov{\Int}(P)$ from the suspension of $\ov{P}$ by a sequence of edge-stellations.

For a bounded poset $P$, let $\Int_{\nonc}(P)$ be the poset of closed intervals $[x,y]$ for which $(x,y)$ is non-contractible, ordered by inclusion.  We note that if $x=y$ or $x\lessdot y$, then $\Delta((x,y))$ is an empty complex, which is non-contractible.

\begin{lemma}\label{lem_interval_noncontractible}
$\ov{\Int}_{\nonc}(P)$ is homotopy equivalent to $\susp(\ov{P})$.
\end{lemma}

\begin{poof}
From the above discussion, it suffices to show that $\ov{\Int}_{\nonc}(P)$ is homotopy equivalent to $\ov{\Int}(P)$.

Let $I_1,\ldots,I_N$ be a linear extension of $\ov{\Int}(P)$.  For $i\geq 0$, let
$$Q_i=\{I_j:\ j\leq i\ \mbox{or}\ I_j\ \mbox{is non-contractible}\}.$$
Then $Q_N=\ov{\Int}(P)$ and $Q_0=\ov{\Int}_{\nonc}(P)$.  If $I_i$ is non-contractible then $Q_{i-1}=Q_i$.

Let $i\geq 0$ and assume $I_i$ is contractible.  Since $I_j\subseteq I_i$ implies $j\leq i$, the subposet $(Q_i)_{<I_i}$ is equal to $\ov{\Int}(I_i)$.  The latter is the suspension of a contractible complex, so it is contractible.  Hence, $Q_{i-1}\simeq Q_i$.  The result now follows by induction.
\end{poof}

\section{Higher Bruhat orders}\label{sec_bruhat}

For $n,d\in\Nbb$, we let $\binom{[n]}{d+1}$ denote the set of $(d+1)$-element subsets of $[n]=\{1,\ldots,n\}$.  A subset $X$ of $\binom{[n]}{d+1}$ is \emph{closed} if for $I\in\binom{[n]}{d-1},\ i,j,k\in[n]\setm I,\ i<j<k$,
$$\mbox{if }I\cup\{i,j\}\in X\mbox{ and }I\cup\{j,k\}\in X\mbox{ then }I\cup\{i,k\}\in X.$$
For $X\subseteq\binom{[n]}{d+1}$, we let $\ov{X}$ be the smallest closed set containing $X$.  If $X$ is a family of subsets of $[n]$ and $P\subseteq[n]$, the \emph{restriction} $X|_P$ of $X$ to $P$ is the subfamily of subsets contained in $P$.  A subset $X$ of $\binom{[n]}{d+1}$ is \emph{consistent} if $X$ and $\binom{[n]}{d+1}\setm X$ are both closed.  Equivalently, $X$ is consistent if for any $(d+2)$-subset $P=\{i_0,\ldots,i_{d+1}\},\ i_0<\cdots<i_{d+1}$,
$$X|_P=\begin{cases}\{P\setm i_{d+1},\ldots,P\setm i_t\}\ &\mbox{for some }t, OR\\\{P\setm i_t,\ldots,P\setm i_0\}\ &\mbox{for some }t\end{cases}.$$

The \emph{higher Bruhat order} $B(n,d)$ is the poset of consistent subsets of $\binom{[n]}{d+1}$ ordered by \emph{single-step inclusion}; that is, $X\leq Y$ if there exists a sequence of consistent subsets $X_0\subseteq\cdots\subseteq X_t$ such that $X=X_0,\ Y=X_t$ and $|X_i\setm X_{i-1}|=1$ for all $i$.  The same set ordered by ordinary inclusion is denoted $B_{\subseteq}(n,d)$.  The posets $B(n,d)$ and $B_{\subseteq}(n,d)$ are both graded with rank function $X\mapsto|X|$ (\cite{ziegler:bruhat} Theorem 4.1(G)).

When $d=1$, $B(n,1)$ is isomorphic to the weak order on the symmetric group, and the two orders $B(n,1)$ and $B_{\subseteq}(n,1)$ coincide.  The weak order $B(n,1)$ is a lattice where the join of $X$ and $Y$ is $\ov{X\cup Y}$.  If $d\geq 2$, the poset $B(n,d)$ may not be a lattice; in particular, $B(6,2)$ is not a lattice (\cite{ziegler:bruhat} Theorem 4.4).  Ziegler proved that $B(n,d)=B_{\subseteq}(n,d)$ when $n-d\leq 4$, but $B(8,3)$ is weaker than $B_{\subseteq}(8,3)$ (\cite{ziegler:bruhat} Theorem 4.5).  In fact, his example in $B(8,3)$ shows that $\ov{X\cup Y}$ need not be the join of $X,Y\in B(n,d)$ even if $\ov{X\cup Y}$ is consistent.

For $X,Y\in B(n,d)$, if $X\subseteq Y$ we define the \emph{ascent set}
$$\Asc(X,Y)=\{I\in Y\setm X:\ X\cup\{I\}\in B(n,d)\}.$$
If $Y=\hat{1}$, we write $\Asc(X)$ for $\Asc(X,Y)$. 


\begin{lemma}\label{lem_ascent_decomp}
Fix $X\in B(n,d)$.  The ascent set $\Asc(X)$ decomposes as the disjoint union $\Asc(X)=A_1\sqcup\cdots\sqcup A_N$ where
\begin{enumerate}
\item\label{lem_ascent_decomp_1} $A_t=\{\{a_{t1}<\cdots<a_{t,d+1}\},\{a_{t2}<\cdots<a_{t,d+2}\},\ldots,\{a_{t,r_t}<\cdots<a_{t,d+r_t}\}\}$ (i.e. $A_t$ is the set of contiguous intervals in the set $\{a_{t1},\ldots,a_{t,r_t+d}\}\subseteq[n]$), and
\item\label{lem_ascent_decomp_2} if $I\in A_s,\ J\in A_t$, $s\neq t$ then $|I\cap J|<d$.
\end{enumerate}
\end{lemma}

\begin{poof}
We first show that any ascent $I\in \Asc(X)$ shares $d$ elements with at most two other ascents of $X$.  Suppose $I,J\in \Asc(X)$ such that $|I\cap J|=d$ with $I<J$ in lexicographic order.  The restriction $X|_{I\cup J}$ is an element of $B(|I\cup J|,d)$ with two ascents, so it must be the bottom element.  Consequently, the $I$ ($J$) is the lex-minimal (lex-maximal) $(d+1)$-subset of $I\cup J$.

Now suppose $J^{\pr}\in \Asc(X),\ J^{\pr}\neq J$ such that $J^{\pr}>I$ in lexicographic order and $|J^{\pr}\cap I|=d$.  Then $J^{\pr}$ is the lex-maximal $(d+1)$-subset of $I\cup J^{\pr}$ by the above argument.  But, $|J\cap J^{\pr}|=d$ and $J,J^{\pr}$ are not at opposite ends of their $(d+1)$-packet, a contradiction.

We have now established that for any $I\in \Asc(X)$, there is at most one $J>I$ in lexicographic order for which $|I\cap J|=d$.  By similar reasoning, there is at most one $L<I$ with $|I\cap L|=d$.  Thus, $\Asc(X)$ decomposes into chains $I_1^t<I_2^t<\cdots<I_{m_t}^t$ where $|I_i^t\cap I_{i+1}^t|=d$ and all other intersections have cardinality strictly less than $d$.

\end{poof}

\begin{lemma}\label{lem_noninv}
If $X\in B(n,d)$, then $X\cap\ov{\Asc(X)}=\emptyset$.
\end{lemma}

\begin{poof}
The ascent set $\Asc(X)$ is a subset of $\binom{[n]}{d}\setm X$, and the latter set is closed.  Hence, $X\cap\ov{\Asc(X)}$ is empty.
\end{poof}

\section{The second higher Bruhat order}\label{sec_bruhat_2}

A \emph{wiring diagram} is a collection of \emph{wires}, continuous piecewise linear curves $C_1,\ldots,C_n$ in $\Rbb^2$, satisfying the following conditions.
\begin{itemize}
\item The projection of $C_i$ onto the first coordinate is bijective.
\item The wires are in order $C_1,\ldots,C_n$ top-to-bottom, sufficiently far to the right.
\item Distinct wires $C_i,C_j$ cross at a unique point.
\item All crossings are transverse.
\end{itemize}

We shall further assume that the wiring diagram is \emph{simple}, meaning there are no common intersections among three or more wires.  In particular, each wire $C_i$ determines a permutation $\pi_i=a_1\cdots a_{n-1}$ of $[n]\setm i$ where if $r<s$ then the first coordinate of $C_i\cap C_{a_r}$ is less than that of $C_i\cap C_{a_s}$.  Two wiring diagrams are considered equivalent if they determine the same sequence of wire permutations $(\pi_i)_{i\in[n]}$.

For $1\leq i<j<k\leq n$, if the crossing of $C_i$ and $C_k$ is below (above) $C_j$, then $\{i,j,k\}$ is an \emph{inversion triple} (\emph{non-inversion triple}).  The map taking a wiring diagram to its set of inversion triples defines a bijection between equivalence classes of simple wiring diagrams with $n$ wires and consistent subsets of $\binom{[n]}{3}$.  A \emph{block} is a set of non-inversion triples of the form $\{\{i_j,i_{j+1},i_{j+2}\}:\ j\in[m]\}$, where $1\leq i_1<\cdots<i_{m+2}\leq n$.


For distinct $i,j,k\in[n]$, the piece $S_{ik}^j$ of $C_j$ between $C_i\cap C_j$ and $C_k\cap C_j$ is called a \emph{segment} of $C_j$.  If $\{i,j,k\}$ is a non-inversion triple, $i<j<k$, then the \emph{floor} of $\{i,j,k\}$ is the segment $S_{ik}^j$.  The \emph{floor of a block} $\Ical$ is the union of the floors of elements of $\Ical$.  A floor is \emph{elementary} if its interior is not intersected by any other wire.  The \emph{height} of a non-inversion triple $\{i,j,k\}$, denoted $\height(\{i,j,k\})$, is the number of wires that pass below the segment $S_{ik}^j$.  If a block $\Ical$ has an elementary floor, then all of its elements have the same height, which we denote $\height(\Ical)$.

\begin{proposition}[\cite{felsner.weil:theorem}]\label{prop_diff_triple}
Let $W$ be a simple wiring diagram with inversion set $X$.  Let $Y\in B(n,2)$ such that $X\subsetneq Y$.
\begin{enumerate}
\item\label{prop_diff_triple_1} There exists an element of $Y-X$ with an elementary floor in $W$. (\cite{felsner.weil:theorem} Lemma 2.2)
\item\label{prop_diff_triple_2} Among those elements of $Y-X$ with an elementary floor, if $I$ is of maximum height, then $X\cup\{I\}$ is consistent.  In particular, $\Asc(X,Y)$ is nonempty.  (\cite{felsner.weil:theorem} Lemma 2.3)
\end{enumerate}
\end{proposition}

\begin{corollary}\label{cor_inclusion}
The second higher Bruhat order $B(n,2)$ is ordered by inclusion; that is, $B(n,2)=B_{\subseteq}(n,2)$ as posets.
\end{corollary}

Given $\Ical\subseteq\binom{[n]}{d+1}$, let $[\Ical]$ denote the union $\bigcup_{I\in \Ical}I$.

\begin{lemma}\label{lem_noninv_wiring}
Let $X\in B(n,2)$ have wiring diagram $W$.  If $\Ical$ is a block with an elementary floor in $W$, then $X\cup\ov{\Ical}$ is not consistent if and only if there exists a wire $p$ intersecting the segments $S^{i_0}_{i_1,i_m}$ and $S^{i_m}_{i_0,i_{m-1}}$ where $[\Ical]=\{i_0<\cdots<i_m\}$.
\end{lemma}

\begin{poof}
Let $\Ical$ be a block with an elementary floor in $W$, and let $p\in[n]-[\Ical]$.  If $X\cup\ov{\Ical}$ is consistent, then for $i\in[\Ical]$ the words $\pi_i$ have the elements of $[\Ical]\setm i$ flipped with the other letters in the same relative order.  Hence, $X\cup\ov{\Ical}$ is consistent if and only if every wire in $[n]-[\Ical]$ does not intersect any segment $S_{ij}^k$ for $i,j,k\in[\Ical]$.

If $X\cup\ov{\Ical}$ is not consistent, then there exists a wire $p$ intersecting some segment $S_{ij}^k$ with $i,j,k\in[\Ical]$.  As $\Ical$ has an elementary floor in $W$, $p$ must intersect the segments $S^{i_0}_{i_1,i_m}$ and $S^{i_m}_{i_0,i_{m-1}}$ by planarity.
\end{poof}

To determine the homotopy type of intervals of $B(n,2)$, we use a stronger version of Proposition \ref{prop_diff_triple}(\ref{prop_diff_triple_2}).

\begin{proposition}\label{prop_diff_block}
Let $W$ be a simple wiring diagram with inversion set $X$.  Let $Y\in B(n,2)$ such that $X\subseteq Y$.  Among the blocks in $Y-X$ with an elementary floor, if $\Ical$ is of maximum height, then $X\cup\ov{\Ical}$ is consistent.
\end{proposition}

\begin{poof}
Let $\Ical$ be a difference block of $Y-X$ of maximum height with an elementary floor, and assume $X\cup\ov{\Ical}$ is not consistent.  Replacing $\Ical$ by a smaller block, we may assume that $X\cup\ov{\Ical^{\pr}}$ is consistent for every block $\Ical^{\pr}$ that is a proper subset of $\Ical$.  Let $[\Ical]=\{i_0,\ldots,i_m\}$ where $i_0<\cdots<i_m$.  By Lemma \ref{lem_noninv_wiring}, there exists a wire $p$ intersecting the segments $S^{i_0}_{i_1,i_m}$ and $S^{i_m}_{i_0,i_{m-1}}$.  By the minimality of $\Ical$, every such wire intersects the subsegments $S^{i_0}_{i_{m-1},i_m}$ and $S^{i_m}_{i_0,i_1}$; see Figure \ref{fig_lem_3_diagram}.

\begin{figure}
\begin{centering}
\begin{tikzpicture}
  \draw (0,0) node[anchor=east]{$i_0$} -- (.5,0) -- (3,5) -- (4,5) -- (5,6) -- (9,6)
        (0,1) node[anchor=east]{$i_1$} -- (.5,1) -- (1,0) -- (2,0) -- (5,3) -- (6.5,3) -- (7,4) -- (9,4)
        (0,2) node[anchor=east]{$i_2$} -- (1,2) -- (1.5,1) -- (2,1) -- (3,0) -- (4,0) -- (6,2) -- (7,2) -- (7.5,3) -- (9,3)
        (0,3) node[anchor=east]{$i_3$} -- (1.5,3) -- (2,2) -- (3,2) -- (5,0) -- (6,0) -- (7,1) -- (7.5,1) -- (8,2) -- (9,2)
        (0,4) node[anchor=east]{$i_{m-1}$} -- (2,4) -- (2.5,3) -- (4,3) -- (7,0) -- (8,0) -- (8.5,1) -- (9,1)
        (0,6) node[anchor=east]{$i_m$} -- (4,6) -- (5,5) -- (6,5) -- (8.5,0) -- (9,0);
  \draw[black!70] (-1,3.5) node[anchor=east]{$p$} -- (1.75,3.5) -- (2.25,2.5) -- (3.5,2.5) -- (4.5,3.5) -- (5,3.5) -- (7,1.5) -- (7.25,1.5) -- (7.75,2.5) -- (9.5,2.5);
  \draw[ultra thick] (.75,.5) -- (2.25,3.5) -- (2.5,3) -- (4,3) -- (6.5,.5) -- (6,0) -- (5,0) -- (4.5,.5) -- (4,0) -- (3,0) -- (2.5,.5) -- (2,0) -- (1,0) -- cycle;
\end{tikzpicture}
\caption{\label{fig_lem_3_diagram}Wire $p$ intersects $S_{i_1i_{m-1}}^{i_0}$ so $X\cup\ov{\Ical\setm\{i_{m-2},i_{m-1},i_m\}}$ is not consistent.}
\end{centering}
\end{figure}
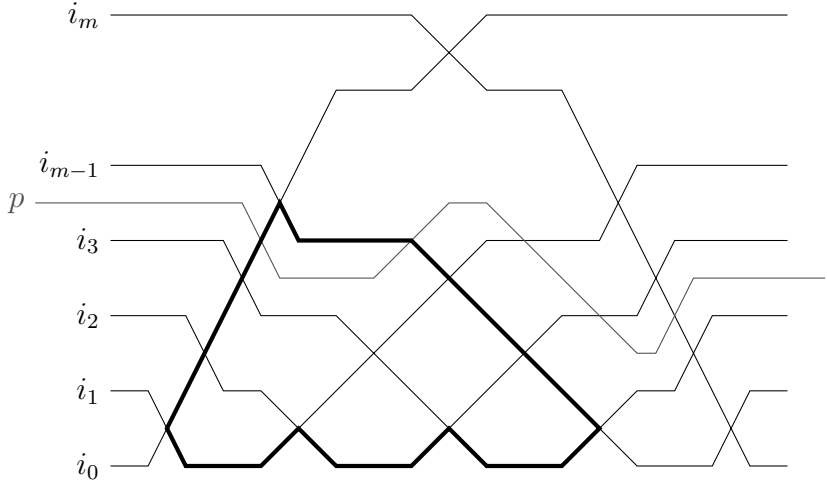

Let $P$ be the set of wires intersecting $S^{i_0}_{i_{m-1},i_m}$ and $S^{i_m}_{i_0,i_1}$.  If $p\in P$, we claim that $i_0<p<i_m$ and $\{i_0,p,i_m\}$ is a difference triple in $Y-X$.  This follows by restriction of $W$ to the wires $\{i_0,p,i_1,i_m\}$.

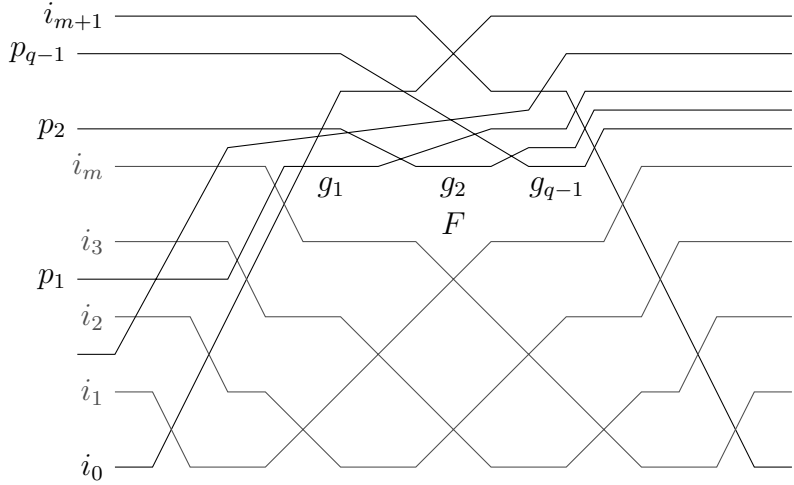
\begin{figure}
\begin{centering}
\begin{tikzpicture}
  \draw (0,0) node[anchor=east]{$i_0$} -- (.5,0) -- (3,5) -- (4,5) -- (5,6) -- (9,6)
        (0,6) node[anchor=east]{$i_{m+1}$} -- (4,6) -- (5,5) -- (6,5) -- (8.5,0) -- (9,0);
  \draw[black!70] (0,1) node[anchor=east]{$i_1$} -- (.5,1) -- (1,0) -- (2,0) -- (5,3) -- (6.5,3) -- (7,4) -- (9,4)
                  (0,2) node[anchor=east]{$i_2$} -- (1,2) -- (1.5,1) -- (2,1) -- (3,0) -- (4,0) -- (6,2) -- (7,2) -- (7.5,3) -- (9,3)
                  (0,3) node[anchor=east]{$i_3$} -- (1.5,3) -- (2,2) -- (3,2) -- (5,0) -- (6,0) -- (7,1) -- (7.5,1) -- (8,2) -- (9,2)
                  (0,4) node[anchor=east]{$i_m$} -- (2,4) -- (2.5,3) -- (4,3) -- (7,0) -- (8,0) -- (8.5,1) -- (9,1);
  \draw (-.5,1.5) node[anchor=east]{} -- (0,1.5) -- (1.5,4.25) -- (5.5,4.75) -- (6,5.5) -- (9,5.5)
        (-.5,2.5) node[anchor=east]{$p_1$} -- (1.5,2.5) -- (2.25,4) -- node[anchor=north]{$g_1$} (3.5,4) -- (5,4.5) -- (6,4.5) -- (6.25,5) -- (9,5)
        (-.5,4.5) node[anchor=east]{$p_2$} -- (3,4.5) -- (3.5,4.25) -- (4,4) -- node[anchor=north]{$g_2$} (5,4) -- (5.5,4.25) -- (6.12,4.25) -- (6.37,4.75) -- (9,4.75)
        (-.5,5.5) node[anchor=east]{$p_{q-1}$} -- (3,5.5) -- (5.5,4) -- node[anchor=north]{$g_{q-1}$} (6.25,4) -- (6.5,4.5) -- (9,4.5);
  \draw (4.5,3.25) node{$F$};
\end{tikzpicture}
\caption{\label{fig_lem_3_gbase}The region $F$, edges $g_1,\ldots,g_{q-1}$, and supporting wires $p_1,\ldots,p_{q-1}$ in the proof of Proposition \ref{prop_diff_block}.}
\end{centering}
\end{figure}

Let $F$ denote the region above the wires $i_1,\cdots,i_{m-1}$, below $i_0,i_{m+1}$ and below all of the wires in $P$.  As shown in Figure \ref{fig_lem_3_gbase}, we label the upper edges of $F$ by $g_0,g_1,\ldots,g_q$ which are supported by the wires $i_0=p_0<p_1<\cdots<p_{q-1}<p_q=i_m$.

We show by induction that one of the $g_j,\ j\in[q-1]$ is the floor of a difference triple in $Y\setminus X$.  We are given that $\{i_0,p_j,i_m\}$ is in $Y$.  Suppose $\{p_{j-1},p_j,i_m\}\in Y$.  Using the packet $\{p_{j-1},p_j,p_{j+1},i_{m+1}\}$ either $\{p_{j-1},p_j,p_{j+1}\}\in Y$ or $\{p_j,p_{j+1},i_m\}\in Y$.  The former case has an elementary floor $g_j$.  Induction on $j$ completes the argument.

Hence there exists a difference triple $\{p_{j-1},p_j,p_{j+1}\}$ with an elementary floor.  But this triple has height strictly greater than that of $\Ical$, a contradiction.
\end{poof}

\begin{lemma}\label{lem_ascent_disjoint}
Let $X\in B(n,2)$ and $\Asc(X)=A_1\sqcup\cdots\sqcup A_N$ as in Lemma \ref{lem_ascent_decomp}.  Suppose $|[A_s]\cap[A_t]|\geq 2$ for some $s\neq t$ and assume $\height(A_s)\leq\height(A_t)$. Then
$$[A_s]\cap[A_t]=\{\min[A_t],\max[A_t]\}.$$
Consequently, $X\cup\ov{A_s}$ is not consistent.
\end{lemma}

\begin{poof}
Let $s,t$ be distinct indices with $|[A_s]\cap[A_t]|\geq 2$ and $\height(A_s)\leq\height(A_t)$.  Let $i,k\in[A_s]\cap[A_t]$ such that $i<k$.  We let $W$ denote a wiring diagram of $X$.

We first show that $\height(A_s)<\height(A_t)$.  By Lemma \ref{lem_ascent_decomp}(\ref{lem_ascent_decomp_2}) there exists $q\in[A_t]\setminus[A_s]$ such that $i<q<k$.  Let $[i,k]\cap[A_s]=\{i=j_0<j_1<\cdots<j_r<j_{r+1}=k\}$ and let $e_{\alpha}$ be the base of $\{j_{\alpha-1},j_{\alpha},j_{\alpha+1}\}$ for $1\leq\alpha\leq r$.  Since $\{j_{\alpha-1},j_{\alpha},j_{\alpha+1}\}$ is an ascent, each $e_{\alpha}$ is a segment.  Let $e=\bigcup_{\alpha} e_{\alpha}$ be the union of these segments.  Then $q$ does not intersect $e_{\alpha}$.  If $q$ is above $e_{\alpha}$, then $\height(A_s)<\height(\{i,q,k\})=\height(A_t)$ as desired.  If $q$ is below $e_{\alpha}$, then $\height(A_s)>\height(\{i,q,l\})=\height(A_t)$, contrary to the hypothesis.

Let $p\in[A_t]\setminus\{i,k\}$.  It remains to show that $i<p<k$.  From this, it follows that $[A_s]\cap[A_t]$ must intersect only at the 2 elements which lie at opposite ends of $[A_t]$.

Suppose to the contrary that $p<i$.  Since $\{p,i,l\}\notin X$, $\pi_p(i)<\pi_p(l)$.  If $\min[A_s]<i$ then the base of $\{p,i,l\}$ includes the base of an ascent $I$ in $A_s$.  By assumption on the height of $A_t$, this implies $I\in A_t$, a contradiction.  If $\min[A_s]=i$ then the base of $\{p,i,l\}$ includes the crossing $i\cap j$ where $j=\min([A_s]\setminus i)$.  Consequently, $\height(A_t)\leq\height(A_s)$, a contradiction.

A symmetric argument shows that $p\ngtr l$, thus completing the proof.
\end{poof}

The following proposition is the key to the proof of Theorem \ref{thm_main}, as described in the introduction.

\begin{proposition}\label{prop_top_atom}
Let $W$ be a simple wiring diagram with inversion set $X$.  Let $Y\in B(n,2)$ such that $X\subsetneq Y$, and let $\Ical\subseteq \Asc(X,Y)$ such that $X\cup\ov{\Ical}$ is consistent.  If $I_0\in \Asc(X,Y)$ is of maximum height in $W$, then $X\cup\ov{\Ical\cup\{I_0\}}$ is consistent.
\end{proposition}


\begin{poof}
We may assume that $I_0\notin\Ical$, as the result is otherwise immediate.

From Lemma \ref{lem_ascent_disjoint}, we know that $\Ical$ uniquely decomposes as a union of blocks $\Ical_1,\ldots,\Ical_m$ such that $|[\Ical_s]\cap[\Ical_t]|\leq 1$ for all $s\neq t$.  By this decomposition, if $s\neq t$, then no packet contains both a subset of $[\Ical_s]$ and of $[\Ical_t]$.   Hence $X\cup\ov{\Ical_s}$ is consistent for all $s$.

Lemma \ref{lem_ascent_decomp} implies that $|I_0\cap[\Ical_s]|=2$ for 0,1, or 2 blocks $\Ical_s$.  We consider each of these cases in turn.

If $|I_0\cap[\Ical_s]|\leq 1$ for all $s$, then $X\cup\ov{\Ical\cup\{I_0\}}=X\cup\ov{\Ical}\cup\{I_0\}$ is consistent.

Suppose $|I_0\cap[\Ical_s]|=2$ for exactly one block $\Ical_s$.  Then $X\cup\ov{\Ical_s\cup\{I_0\}}$ is consistent by Proposition \ref{prop_diff_block}.  By Lemma \ref{lem_ascent_disjoint}, we deduce that $|[\Ical_s\cup\{I_0\}]\cap[\Ical_t]|\leq 1$ for $t\neq s$, so $X\cup\ov{\Ical\cup\{I_0\}}$ is consistent in this case.

Finally, assume that $|I_0\cap[\Ical_s]|=2$ and $|I_0\cap[\Ical_t]|=2$ for two blocks $\Ical_s,\Ical_t$.  Then $X\cup\ov{\Ical_s\cup\{I_0\}\cup\Ical_t}$ is consistent by Proposition \ref{prop_diff_block}.  Lemma \ref{lem_ascent_disjoint} implies that $|[\Ical_s\cup\{I_0\}\cup\Ical_t]\cap[\Ical_u]|\leq 1$ if $u\neq s$ and $u\neq t$.  Therefore, $X\cup\ov{\Ical\cup\{I_0\}}$ is consistent.

\end{poof}

\section{Proof of Theorem \ref{thm_main}}\label{sec_main}

Let $V,W$ be wiring diagrams such that $W$ is simple and $V$ is non-simple.  We say $W$ is \emph{incident} to $V$ if $V$ may be obtained by moving the wires of $W$ to a more special position.  More precisely, $W$ is incident to $V$ if the associated oriented matroid of $V$ is a weak map image of the oriented matroid associated to $W$; see \cite[Section 7.7]{bjorner.lasVergnas.ea:oriented} for background on weak maps.  An interval $(X,Y)$ of $B(n,2)$ is called \emph{facial} if the closed interval $[X,Y]$ is the set of inversion sets of simple wiring diagrams incident to some fixed wiring diagram.

The following lemma follows from the proof of Lemma \ref{lem_ascent_disjoint}.

\begin{lemma}\label{lem_facial_intervals}
Let $X,Y\in B(n,2)$ such that $X<Y$.  Then $(X,Y)$ is facial if and only if $Y=X\cup\ov{\Asc(X,Y)}$ and $X\cup\ov{\Ical}\in B(n,2)$ for $\Ical\subseteq\Asc(X,Y)$.
\end{lemma}

\begin{theorem}\label{thm_main_detailed}
Let $X,Y\in B(n,2)$ such that $X<Y$.  If $[X,Y]$ is facial, then $(X,Y)$ is homotopy equivalent to a sphere of dimension $(|\Asc(X,Y)|-2)$.  Otherwise, $(X,Y)$ is contractible.
\end{theorem}

\begin{poof}
Assume the statement holds for intervals $(X,Z)$ with $X\leq Z<Y$.  By Lemma \ref{lem_noncontractible},
$$(X,Y)\simeq (X,Y)_{nonc}$$
where $(X,Y)_{nonc}=\{Z\in(X,Y)\ |\ (X,Z)$ is not contractible$\}$.  Since $B(n,2)$ is ordered by inclusion, $X\cup\ov{\Ical}\leq Y$ whenever $\Ical$ is a subset of $\Asc(X,Y)$ such that $X\cup\ov{\Ical}$ is consistent.  By the inductive hypothesis,
$$(X,Y)_{nonc}=\{Z\in(X,Y)\ |\ (X,Z)\ \mbox{is facial}\}.$$

Suppose $(X,Y)$ is facial.  By Lemma \ref{lem_facial_intervals}, $Z\in(X,Y)_{nonc}$ if and only if $Z=X\cup\ov{\Ical}$ for some non-empty proper subset $\Ical$ of $\Asc(X,Y)$.  Hence $(X,Y)_{nonc}$ is the face poset of the boundary of a simplex.  Thus, $(X,Y)$ is homotopy equivalent to a sphere of dimension $|\Asc(X,Y)|-2$.

Now assume that $(X,Y)$ is not facial.  By Lemma \ref{lem_facial_intervals}, $(X,Y)_{nonc}$ is the face poset of a simplicial complex over $\Asc(X,Y)$.  By Proposition \ref{prop_top_atom}, this simplicial complex has a cone point $I_0\in \Asc(X,Y)$, so $(X,Y)_{nonc}$ is contractible.
\end{poof}

We let $\omega(n,2)$ denote the set of all wiring diagrams on $n$ wires up to equivalence. This forms a poset under the incidence relation. Alternatively, $\omega(n,2)$ may be viewed as the poset of facial intervals of $B(n,2)$, ordered by inclusion. By Lemma \ref{lem_interval_noncontractible} and Theorem \ref{thm_main_detailed}, we deduce that $\ov{\omega}(n,2)$ is homotopy equivalent to $\Sbb^{n-3}$. As observed in the introduction, the homotopy type of $\ov{\omega}(n,2)$ was already computed in \cite{sturmfels.ziegler:extension} by other means.


\section*{Acknowledgements}

The author thanks Vic Reiner for suggesting this problem and providing many references.

\small
\bibliographystyle{plain}
\bibliography{bib_higher_bruhat}

\end{document}